\date{April 17, 2013}
\documentclass[12pt]{amsart}
\usepackage{latexsym,amsmath,amsfonts,amscd,amssymb}
\usepackage{graphics}
\newtheorem{theorem}{Theorem} %[subsection]
\newtheorem{lemma}[theorem]{Lemma}
\newtheorem{proposition}[theorem]{Proposition}
\newtheorem{corollary}[theorem]{Corollary}

\newtheorem{definition}[theorem]{Definition}

\theoremstyle{remark}
\newtheorem{remark}[theorem]{Remark}

\newcommand{\unit}{\textbf{1}}

\newcommand{\supp}{\operatorname{supp}}

\newcommand{\ord}{\operatorname{ord}}

\newcommand{\cC}{{\mathcal C}}

\newcommand{\cL}{{\mathcal L}}

\newcommand{\CC}{{\mathbb C}}

\newcommand{\NN}{{\mathbb N}}

\newcommand{\QQ}{{\mathbb Q}}
\newcommand{\RR}{{\mathbb R}}

\newcommand{\ZZ}{{\mathbb Z}}

\newcommand{\bk}{{\mathbf{k}}}

\title{On the genus of meromorphic functions}


\author[V. Mu\~{n}oz]{Vicente Mu\~{n}oz}
\address{Facultad de
Matem\'aticas, Universidad Complutense de Madrid, Plaza de Ciencias
3, 28040 Madrid, Spain}
\email{vicente.munoz@mat.ucm.es}

\author[R. P\'{e}rez Marco]{Ricardo P\'{e}rez Marco}
\address{CNRS, LAGA UMR 7539, Universit\'e Paris XIII,
99, Avenue J.-B. Cl\'ement, 93430-Villetaneuse, France}
\normalsize\email{ricardo@math.univ-paris13.fr}

\thanks{Partially supported through Spanish MICINN grant MTM2010-17389.}

\begin{document}

\maketitle

\begin{abstract}
We define the class of Left Located Divisor (LLD) meromorphic functions and their vertical order $m_0(f)$ and 
their convergence exponent $d(f)$. When $m_0(f)\leq d(f)$ we prove that their Weierstrass genus is minimal. This 
explains the phenomena that many classical functions have minimal Weierstrass genus, for example Dirichlet series, 
the $\Gamma$-function, and trigonometric functions.
\end{abstract}

\subsection{LLD meromorphic functions.}

Meromorphic functions $f$ on $\CC$, of the variable $s\in \CC$, considered in this article are assumed to be of finite 
order $o=o(f)$. We recall that the order $o(f)$ is defined as
$$
o(f)=\limsup_{R\to +\infty} \frac{\log \log ||f||_{C^0(B(0,R))}}{\log R}  \ .
$$

We study in this article Dirichlet series, and more generally 
the class of meromorphic of finite order with Left Located Divisor (LLD), which we call LLD meromorphic functions:

\begin{definition} {\textbf{(LLD meromorphic functions)}}
A LLD meromorphic function is a function $f$ of finite order and  left located divisor 
$$
\sigma_1 =\sup_{\rho \in f^{-1}(\{0,\infty\})} \Re \rho < +\infty \ .
$$
\end{definition}

The properties that we establish in this article are invariant by a real translation. Thus 
considering $g(s)=f(s+\sigma_1)$ instead of $f$ we will assume that $\sigma_1=0$.

Examples of LLD meromorphic functions are Dirichlet series, that we normalize in this article 
such that $f(s)\to 1$ when $\Re s\to +\infty$. A Dirichlet series is of the form
\begin{equation}\label{eqn:1}
 f(s)=1+\sum_{n\geq 1} a_n \ e^{-\lambda_n s} \ ,
\end{equation}
with $a_n \in \CC$ and
$$
0< \lambda_1 < \lambda_2 < \ldots
$$
with $(\lambda_n)$ a discrete set, that is either finite  
or $\lambda_n \to +\infty$, and such that we have a half
plane of absolute convergence, i.e., for some $\bar \sigma \in \RR$ we have
$$
\sum_{n\geq 1} |a_n | \ e^{-\lambda_n \bar \sigma} <+\infty \, .
$$
We refer to \cite{HR} for classical background on Dirichlet series.

\subsection{Convergence exponent.}

We denote by $(\rho)$ the set of zeros and poles
of $f$, and the integer $n_\rho$ is the multiplicity of $\rho$ (positive for zeros and
negative for poles, with the convention $n_\rho =0$ if $\rho$ is neither a zero nor a pole). 

\begin{definition} \textbf{(Convergence exponent)}
 The convergence exponent of $f$ is the minimum
integer $d=d(f)\geq 0$ such that
$$
\sum_{\rho \not= 0} |n_\rho| \, |\rho|^{-d} < +\infty\ .
$$
\end{definition}

We have $d=0$ if and only if $f$ has a finite divisor, i.e. it 
is a rational function multiplied by the exponential of a polynomial, otherwise
$d\geq 1$. 

It is classical that the convergence exponent satisfies $d \leq [o]+1$ (see \cite{A}), thus it is finite 
for functions of finite order, but there is no upper bound 
of the order by the convergence exponent since we can always multiply by $\exp P$, where $P$ is a polynomial, 
increasing the order without changing the divisor, hence keeping the same convergence exponent.

\subsection{Genus.}

When $f$ is a meromorphic function of finite order we have the Hadamard factorization of $f$ (see \cite{A}, p.208)
$$
f(s)=s^{n_0} e^{Q_f(s)} \prod_{\rho \not=0 } \left ( E_{d-1} (s/\rho )\right )^{n_\rho} \ ,
$$
where 
$$
E_n(z)=(1-z) e^{z+\frac12 z^2 +\ldots +\frac1n z^n} \ ,
$$
and $Q_f$ is a polynomial, the Weierstrass polynomial, 
uniquely defined up to the addition of an integer multiple of $2\pi i$. 

The \textit{discrepancy polynomial} of the meromorphic function $f$ is
$$
P_f=-Q_f' \ .
$$ 

We define the Hadamard part of $f$ as 
\begin{equation}\label{eqn:H}
f_H(s)=s^{n_0} \prod_{\rho \not=0 } \left ( E_{d-1} (s/\rho )\right )^{n_\rho} \ .
\end{equation}
Note that $f'/f= f_H'/f_H - P_f$.

The degree $g_W=\deg Q_f$ is the Weierstrass genus. The genus 
of $f$ is defined as the integer
$$
g=g(f)=\max (g_W(f) , g_H(f)) \ .
$$
where $g_H(f)=d(f)-1$ is the Hadamard genus, which is the degree of the polynomials in the exponential 
of the  factors $E_n(z)$.
From the definition we have $d\leq g+1$, and 
$g\leq o \leq g+1$ (see \cite{A}, p.209).

We set the following useful definition:

\begin{definition}\textbf{(Hadamard and Weierstrass type)}
 A meromorphic function $f$ is of Hadamard type when $g(f)=g_H(f)=d(f)-1\geq g_W(f)$. It is of Weierstrass type when 
$g(f)=g_W(f)>g_H(f)$. 
\end{definition}
 
Many classical functions are of Hadamard type. One of the purposes of the article is to explain why this holds.

\subsection{Vertical order.}

For a LLD meromorphic function we look at the growth of its logarithmic
derivative on the right half plane. This growth is always polynomial (proof in Appendix 1). 

\begin{proposition} \label{prop:estimate}
The logarithmic derivative of a LLD meromorphic function has polynomial 
growth on a right half plane, i.e. for 
$\sigma_2 >\max (0, \sigma_1)$, and for $\Re s > \sigma_2$,
$$
\left |\frac{f'(s)}{f(s)} \right | \leq C_0 |s|^{\max (d, g_W-1)}  \ ,
$$
more precisely we have
$$
\left |\frac{f'_H(s)}{f_H(s)} \right | \leq C_0 |s|^{d} 
$$
\end{proposition}

\begin{remark}
The exponent $d$ is best possible in the last estimate (see the example constructed in Appendix 2).
\end{remark}

We define the {\it vertical order} as follows:

\begin{definition}\textbf{(Vertical order)}
The vertical order of a meromorphic function $f$ with left located divisor is the 
minimal integer $m_0=m_0(f) \geq 0$
such that for $c> \sigma_1$, $c \not= 0$,
$$
|c+it|^{-m_0} \frac{f'}{f} (c+it) \in L^1(\RR ) \ . 
$$ 
\end{definition}

\begin{lemma} \label{lem:ccc}
This definition does not depend on the choice of $c$.
\end{lemma}

This Lemma is proved in Appendix 3.
 
From the estimate in Proposition \ref{prop:estimate}
we have that $m_0(f_H) \leq d+2$. But we can do better:

\begin{proposition}\label{prop:vertical}
We have $m_0(f_H) \leq d+1$.
\end{proposition}

For a Dirichlet series normalized as in (\ref{eqn:1}) we have that $f(s)\to 1$ and 
$f'(s)\sim -\lambda_1 a_1 e^{-\lambda_1 s}$ uniformly with 
 $\Re s \to +\infty$,
thus $m_0(f) =2$. 

In this article we say that 
a distribution has order $n$ if $n$ is the minimal integer such that it 
is the $n$-th derivative of a continuous function (there is no consensus in 
the classical literature on the definition of order of a distribution, for example see \cite{S} and \cite{Z}).
Proposition \ref{prop:estimate} implies that the inverse 
Laplace transform $\cL^{-1}(f'/f)$ is a distribution of finite order. 
This is because we have an explicit formula for the inverse Laplace transform. 
We recall (see \cite{Z}) that 
$$
\cL^{-1}(F) (t) =\frac1{2\pi} \int_\RR F(c+iu) e^{(c+iu)t} \ du \ ,
$$
if the integral is convergent, and 
$$
\cL^{-1}(F) (t) =\cL_c^{-1}(F) (t)=\frac1{2\pi} \frac{D^n}{Dt^n}\int_\RR \frac{F(c+iu)}{(c+iu)^n} e^{(c+iu)t} \ du \ ,
$$
in general (the derivative is taken in distributional sense) which holds for some $n$ when 
$F$ is holomorphic with polynomial growth on $\{ \Re s > \sigma_2\}$ and it is independent of 
$c >\sigma_2> \sigma_1$.

A closely related integer to the vertical order is the \textit{distributional vertical order}. 

\begin{definition}\textbf{(Distributional vertical order)}
The distributional vertical order of a LLD meromorphic function $f$ is the 
minimal integer $m \geq 0$ such that the inverse Laplace transform
$$
\cL^{-1}(f'/f)
$$
is a distribution of order $m$.
\end{definition}

It is clear that:

\begin{proposition}
 We have
$m(f)\leq m_0(f)$.
\end{proposition}

\subsection{Main results.}

\begin{theorem} \label{thm:main}
For a LLD meromorphic function $f$ we have that if 
$m(f)\not= g_W(f)+1$ then $f$ is of Hadamard type, i.e. $g_W(f)\leq g_H(f)=g(f)$.

Moreover, any  Dirichlet series $f$ is of Hadamard type, i.e. $g_W(f)\leq g_H(f)=g(f)$ unconditionally.
\end{theorem}

\begin{corollary} \label{cor:cor1}
 If a LLD meromorphic function $f$ is of Weierstrass type then $m(f)=g_W(f)+1$.
\end{corollary}

The next corollary gives an analytic criterium to determine if a meromorphic function is of Hadamard type.

\begin{corollary} \label{cor:cor2}
 If $m_0(f) \leq d(f)$ then $f$ is of Hadamard type.
\end{corollary}

The same argument used in the proof of the main theorem gives:

\begin{theorem}\label{thm:thm2}
 Let $f$ be a non-constant Dirichlet series. Then we have
$$
d(f)\geq 2 \ ,
$$
and
$$
o(f)\geq 1 \ .
$$
\end{theorem}

Before proving these results we need to introduce the Newton-Cramer distribution and Poisson-Newton formula.

\subsection{Newton-Cramer Distribution.}

In \cite{MPM} we associate to the divisor $\mathrm{div}(f)= \sum n_\rho\,\rho$ its Newton-Cramer 
distribution, which is given by the series
$$
W(f)= \sum n_\rho \, e^{\rho t} \ 
$$ 
on $\RR_+^*$. This sum is only converging in $\RR_+^*$ in the distribution sense. 
The distribution $W(f)$ vanishes in $\RR_-^*$, 
and has some structure at $0$. The precise definition follows 
(we assume, in order to simplify, that $\rho =0$ is not part of the divisor).

\begin{definition}\textbf{(Newton-Cramer distribution)}
The Newton-Cramer distribution is
$$ 
W(f)=\frac{D^d}{Dt^d} \left ( L_d(t)         \right ),
$$
where $L_d$ is the continuous function on $\RR$ defined on $\RR_+$ by
$$
L_d(t) =\sum_{\rho\not= 0} \frac{n_\rho}{\rho^d} (e^{\rho t}-1)\unit_{\RR_+} \ .
$$
\end{definition}

It is easy to see that the sum converges for $t\geq 0$.

In this article, only the order of distributions plays a role, and the space of test funcitons 
for which the distribution belong to the dual is not so important. 
The distribution $W(f)$ is Laplace transformable, that is, it can be paired with $e^{-st}$ on $\RR_+$, 
on some half-plane $\Re s>\sigma_0$. Hence, the appropriate space of distributions to
use is the dual of the space of $\cC^\infty$ functions on $\RR$ which decay
faster than $C e^{\alpha |t|}$, for some $C>0$, $\alpha>0$.

The main property of the Newton-Cramer distribution that we need follows from its definition:

\begin{proposition}
 The Newton-Cramer distribution is the $d$-th derivative of a continuous function.
\end{proposition}

\subsection{Poisson-Newton formula.}

The Newton-Cramer distribution of $f$ is linked to the inverse Laplace transform of the logarithmic derivative $f'/f$
by the Poisson-Newton formula (see \cite{MPM}):

\begin{theorem} \label{thm:PN} \textbf{(Poisson-Newton formula)}
For a LLD meromorphic function $f$ we have on $\RR$
$$
W(f)=\sum_{l=0}^{g_W-1} c_l \delta_0^{(l)} + \cL^{-1}(f'/f) \ ,
$$ 
where $P_f(s)=c_0+c_1 s+ \ldots +c_{g_W-1} s^{g_W-1}=-Q_f'(s)$ is the discrepancy polynomial.
\end{theorem}

When $f$ is a Dirichlet function, the Laplace transform $\cL^{-1}(f'/f)$ is purely atomic with atoms 
in $\RR_+^*$. We can compute it explicitely as follows. On the half plane $\Re s >  \sigma_1$, $\log f(s)$ is well
defined taking the principal branch of the logarithm. Then we can
define the coefficients $(b_{\bk})$ by 
\begin{equation} \label{eqn:bn}
-\log f(s)=-\log \left ( 1+ \sum_{n\geq 1} a_n \ e^{-\lambda_n s}\right )
=\sum_{\bk \in \Lambda} b_{\bk} \, e^{-\langle \boldsymbol{\lambda} , \bk \rangle s}
 \ ,
 \end{equation}
where $\Lambda=\{ \bk=(k_n)_{n\geq 1} \, | \, k_n \in \NN, ||\bk||=\sum | k_n |<\infty, ||\bk|| \geq 1\}$,
and
$\langle \boldsymbol{\lambda} , \bk \rangle = \lambda_1k_1+\ldots + \lambda_{l}k_{l}$, where
$k_n=0$ for $n>l$.
Note that the coefficients $(b_{\bk})$ are polynomials on the $(a_n)$. More precisely, we have
\begin{equation} \label{eqn:bs}
 b_\bk= \frac{(-1)^{||\bk||}}{||\bk||} \, \frac{||\bk|| !}{\prod_j k_j!}\, \prod_j a_j^{k_j}\, .
\end{equation}

Note that if the $\lambda_n$ are $\QQ$-dependent then there are repetitions in
the exponents of (\ref{eqn:bn}).

Since $\cL (e^{-\lambda s})=\delta_{\lambda}$, we have
$$
\cL^{-1}(f'/f)= \sum_{\bk \in \Lambda} \langle \lambda , \bk
\rangle \, b_{\bk} \ \delta_{\langle \boldsymbol{\lambda} ,\bk\rangle } \ .
$$

Note in particular that $\supp \cL^{-1}(f'/f) \subset [\epsilon, +\infty[$ for some $\epsilon >0$.

\subsection{Proof of the main results.}

The proof of Theorem \ref{thm:main} 
consists on inspecting the orders of the distributions in both sides of the Poisson-Newton equation:
$$
W(f)=\sum_{l=0}^{g_W-1} c_l \delta_0^{(l)} + \cL^{-1}(f'/f) \ .
$$ 

We will use that for two distributions $U$ and $V$, if $\ord (U)\not= \ord (V)$
then 
$$
\ord (U+V) =\max (\ord (U) , \ord (V) ) \ .
$$

The left hand side is of order $\leq d$ since $W(f)$ is the $d$-th derivative of a continuous function. 

Observe that the Dirac $\delta_0$ is of order $2$, and $\delta_0^{(l)}$ 
is of order $l+2$. In particular, the first term of the right hand side in Poisson-Newton equation is of 
order $g_W+1$. 

The second term of the right hand side is of order $m(f)$ by definition of $m(f)$. 

To prove 
Theorem \ref{thm:main} we assume first that $m < g_W+1$. Then the order of the right hand side 
in Poisson-Newton formula is $g_W +1$. Therefore $d\geq g_W+1$ so $g=g_H\geq g_W$ and $f$ is of Hadamard type. 

We look at the second case when $m> g_W+1$. Then the order of the right hand side is $m$, thus comparing with the 
left hand side, we get $d\geq m > g_W+1$, therefore $g=g_H > g_W$ and $f$ is again of Hadamard 
type. This proves the first statement of the main theorem.

For a Dirichlet series $f$ the distribution $\cL^{-1}(f'/f)$ has support away from $0$, therefore looking at the local 
order at $0$ (which is smaller or equal than the global order) 
of both sides of the equation we get that $d\geq g_W+1$ unconditionally. 
This gives $g=g_H \geq g_W$ and $f$ 
is always of Hadamard type. This ends the proof of Theorem \ref{thm:main}.

Now Corollary \ref{cor:cor1} is a direct application of the main theorem.

For Corollary \ref{cor:cor2} we observe that $m_0(f) \leq d(f)$ gives $m(f)\leq m_0(f) \leq d (f)\leq g(f) +1$. If the last 
inequality is an equality, then $f$ is of Hadamard type and we are done. Otherwise we have  
$g=g_W$ and $m(f) < g_W(f)+1$ 
and using the main theorem 
we get also that $f$ is of Hadamard type, and $g=g_W=g_H=d-1$.

For the proof of Theorem \ref{thm:thm2},  we inspect as before the order of the distributions in the Poisson-Newton-formula. The right hand side contains Dirac
distributions at the frequencies, hence it is at least a second derivative of a continuous function. In the left 
hand side we have $W(f)$ that is the $d$-th derivative of a continuous function. This gives $d\geq 2$.

Also we know that $d\leq o +1$, hence $o\geq 1$.

\subsection{Proof of Poisson-Newton formula.}

Let us prove Theorem \ref{thm:PN}.
We start from the Hadamard factorization of $f$ (assuming that $\rho=0$ is not 
part of the divisor in order to simplify).
$$
f(s)= e^{Q_f(s)} \prod_{\rho  } \left ( E_{d-1} (s/\rho )\right )^{n_\rho} \ ,
$$
We take its logarithmic derivative:
\begin{align} \label{eqn:aaa}
f'/f &= -P_f+\sum_{\rho} n_\rho \frac{E'_{d-1}(s/\rho)}{E_{d-1}(s/\rho)}\notag\\
&=-P_f+\sum_{\rho}  n_\rho \left (\frac{1}{\rho-s} + 
\sum_{l=0}^{d-2}\frac{s^l}{\rho^{l+1}} \right )
\end{align}

Since for $l\geq 0$ 
$$
\cL(\delta_0^{(l)}) =s^{l} \ ,
$$
the polynomial $P_f$ is the Laplace transform
$$
P_f=\cL \left (c_0 \delta_0 +c_1 \delta'_0 +\ldots +c_{g-1} \delta_0^{(g-1)}\right )  \ .
$$
It remains to prove that 
$$
\cL (W(f))= \sum_{\rho}  n_\rho \left (\frac{1}{\rho-s} + \sum_{l=0}^{d-2}\frac{s^l}{\rho^{l+1}} \right )\ .
$$
We have
$$
\frac{D^d}{Dt^d} \left ( (e^{\rho t}-1) \unit_{\RR_+} \right ) = \rho^d e^{\rho t} \unit_{\RR_+} +
\sum_{l=0}^{d} \rho^{d-1-l} \delta_0^{(l)} \ ,
$$
thus for a finite set $A$ of zeros and poles of the divisor, we have
\begin{align*}
W_A(f) &=\sum_{\rho \in A} n_\rho \rho^{-d} 
\frac{D^d}{Dt^d} \left ( e^{\rho t}-1 \right ) \unit_{\RR_+}  \\  
&=\sum_{\rho \in A} n_\rho \left ( e^{\rho t} \unit_{\RR_+}  +
\sum_{l=0}^{d} \rho^{-1-l} \delta_0^{(l)} \right ) \ .
\end{align*}

Now we have
$$
\cL(e^{\rho t} \unit_{\RR_+}) =\frac{1}{\rho-s} \ ,
$$
so
$$
\cL \left ( W_A(f) \right ) =\sum_{\rho \in A} n_\rho \left (\frac{1}{\rho-s}  + 
\sum_{l=0}^{d-2}\frac{s^l}{\rho^{l+1}} \right )  \ ,
$$
and we are done taking the inverse Laplace transform.

\subsection{Application to trigonometric functions.}

We check that the sine function is of Hadamard type. For this it is enough 
to consider the hyperbolic sine function which is an entire function of order $1$,
$$
f(s)=\sinh (s) =\frac{e^s-e^{-s}}{2i} \ . 
$$
The zeros are, for $k\in \ZZ$,
$$
\rho_k = \pi i k \ ,
$$
thus $\sinh$ is a LLD entire function and $d(f)=2$. Also we have
$$
f'(s)/f(s)=\cosh (s)/\sinh(s)=\frac{1+e^{-2s}}{1-e^{-2s}} \to 1
$$
when $\Re s \to +\infty$. Therefore $m_0(f)=2$. 

Using Corollary \ref{cor:cor2} we get

\begin{proposition}
 The function $f(s)=\sinh(s)$ is of Hadamard type.
\end{proposition}

This is something that we know from its Hadamard factorisation (due to Euler)
$$
\sinh (s)= s \prod_{k\in \ZZ^*} \left (1-\frac{s}{\pi ik}\right ) e^{\frac{s}{\pi ik}} \ .
$$

\begin{corollary}
 The function $f(s)=\sin(s)$ is of Hadamard type.
\end{corollary}

\subsection{Application to the $\Gamma$ function.}

We check, without computing its Hadamard factorisation, that the classical $\Gamma$ function is 
of Hadamard type.

The $\Gamma$-function has no zeros and has simple poles at the negative integers. Thus it is a LLD 
meromorphic function and $d=2$. Stirling formula indicates that we must have $m_0(\Gamma) =2$ and 
we check:

\begin{lemma}
 For $c>0$ we have for some constant $C_0>0$ and for $|u|\geq 1$
$$
\left |\frac{\Gamma'}{\Gamma } (c+iu)\right | \leq  \log |u| +C_0 \,
$$
and $m_0(\Gamma)=2$.
\end{lemma}

The classical Stirling's asymptotics holds in a right cone, but we need the estimate in a vertical line,
thus we need to refine the classical estimate. We start with Binet's second formula (see \cite{WW} p.251):
$$
\log \Gamma (s) = \left (s-\frac12 \right ) \log s -s +\frac12 \log (2\pi ) +\varphi (s) \ ,
$$
where 
$$
\varphi(s)=2 \int_0^{+\infty } \frac{\arctan (t/s)}{e^{2\pi t}-1} \, dt \ .
$$

Taking one derivative in the above formula, we get an identity for the digamma function
$$
\psi(s)=\frac{\Gamma'(s)}{\Gamma(s)}=\log s -\frac{1}{2s}+\varphi'(s) \ ,
$$
and 
$$
\varphi'(s)= -2\int_0^{+\infty} \left (\frac{s}{s^2+t^2} \right )\left (\frac{t}{e^{2\pi t}-1}\right ) \, dt \ .
$$
Since
$$
\int_0^{+\infty} \frac{t}{e^{2\pi t}-1} \, dt =\frac{B_2}{4}=\frac{1}{24} \ ,
$$
and if $s=c+iu$ with $c=\Re s >0$,
$$
\left| \frac{s}{s^2+t^2} \right| \leq \frac{1}{|c|} \ ,
$$
we have the estimate
$$
|\varphi'(s)| \leq \frac{1}{24 |c|} \ ,
$$
so $|\psi(s)| \leq \log |s|+C_0$, and the lemma follows.

\medskip

Now we have $m_0(\Gamma) =2 \leq d(\Gamma)=2$ so the application of Corollary \ref{cor:cor2} gives:

\begin{proposition}
 The $\Gamma$ function is a meromorphic function of Hadamard type.
\end{proposition}

\subsection{Application to the Riemann zeta function.}

The Riemann zeta function is a Dirichlet series,
$$
\zeta(s)=\sum_{n=1}^{+\infty } n^{-s} \ ,
$$
and has a meromorphic extension of order 1 to the whole complex plane.
So it is a LLD meromorphic function.

We that $d(\zeta)\leq 2$ from the order $1$, and $d(\zeta)\geq 2$ for the summation of the trivial 
zeros that lie at the even negative integers, thus $ d(\zeta)= 2$.

The logarithmic derivative is bounded on vertical lines and so $m_0(\zeta)=2$. Again, using Corollary \ref{cor:cor2}
we get:

\begin{proposition}
 The Riemann zeta function $\zeta$ is a meromorphic function of Hadamard type.
\end{proposition}

\subsection{Appendix 1: Proof of propositions \ref{prop:estimate} and  \ref{prop:vertical}.}  

We start by considering the analogue of (\ref{eqn:aaa}) centered at $\sigma_1$. This is
 $$
 f'/f =-P_f+\sum_{\rho}  n_\rho \left (\frac{1}{\rho-s} + 
\sum_{l=0}^{d-2}\frac{(s-\sigma_1)^l}{(\rho-\sigma_1)^{l+1}} \right )
 $$
We write $f'/f=-P_f + G$, where 
 $G(s)=\sum n_\rho \, g_\rho(s)$, where 
 $$
g_\rho(s)=
\frac{1}{\rho-s} + 
\sum_{l=0}^{d-2}\frac{(s-\sigma_1)^l}{(\rho-\sigma_1)^{l+1}}=
\frac{(s-\sigma_1)^{d-1}}{(\rho-\sigma_1)^{d-1}} \frac{1}{\rho-s}
 $$
In order to prove Proposition \ref{prop:estimate}, we need to bound 
$|g_\rho(s)|\leq C|s-\sigma_1|^d |\rho-\sigma_1|^{-d}$, for a uniform 
constant $C$, since $\sum n_\rho\, |\rho-\sigma_1|^{-d}<\infty$. For
this we need to bound uniformly
 $$
g_\rho(s)=\frac{\rho-\sigma_1}{(s-\sigma_1)(\rho-s)}
 $$
on the half-plane $\Re s>\sigma_2$. 

If $|\sigma_1-s|\leq \frac12 |\rho-\sigma_1|$ then $|\rho-s|\geq |\rho-\sigma_1|-
|\sigma_1-s| \geq \frac12 |\rho-\sigma_1|$. So $|g_\rho(s)|\leq \frac{2}{|s-\sigma_1|}\leq
C$, as $|s-\sigma_1|\geq \sigma_2-\sigma_1$ is bounded below.

If $|\sigma_1-s|\geq \frac12 |\rho-\sigma_1|$ then $|g_\rho(s)|\leq \frac{2}{|\rho-s|}\leq C$,
as $|s-\rho|\geq \sigma_2-\sigma_1$ is bounded below.

\bigskip

We prove now Proposition \ref{prop:vertical}.
Fix $c>\sigma_1$, and let $a=c-\sigma_1>0$.
We need to see that $G(s) |s-\sigma_1|^{-d-1}$ is integrable, and it is enough to see
that
 \begin{equation} \label{eqn:cota}
\int_{L_c} \frac{|\rho- \sigma_1|}{|s-\sigma_1|^2 |\rho-s|} ds
   \end{equation}
is bounded uniformly on $\rho$, for $L_c=c+i\RR$.

We consider two sets:
\begin{itemize}
\item  $A=\{s\in L_c\, | \, |\rho-\sigma_1|\leq \frac32 |\rho-s|\}$. This is an infinite portion of
$L_c$. The integral is bounded by 
 $$
 \frac32 \int_{L_c} \frac{1}{|s-\sigma_1|^2} ds<\infty.
 $$
\item $B=\{s\in L_c\, | \,|\rho-\sigma_1|\geq \frac32|\rho-s|\}$. This is the intersection of a disc of radius 
$\frac23 |\rho-\sigma_1|$ with $L_c$. So its length is bounded by $\frac43 |\rho-\sigma_1|$. 
The integral there is bounded by 
 $$
 \frac43 \max\left\{ \frac{|\rho- \sigma_1|^2}{|s-\sigma_1|^2 |\rho-s|} \, |\, s\in B \right\}.
 $$
We have that $|\rho-s|\geq a$, so $|\rho-s|^{-1/2}\leq \frac{1}{\sqrt{a}}\leq \frac12$, for $a\geq 4$.
Then $|\rho-s|+ |\rho-s|^{1/2}\leq \frac32|\rho-s|$ and
 $$
|\rho-s|+ |\rho-s|^{1/2}\leq 
|\rho-\sigma_1| \leq |\rho-s|+|s-\sigma_1|.
 $$ 
So $|\rho-s|^{1/2}\leq |s-\sigma_1|$ and
 $$
\frac{|\rho- \sigma_1|^2}{|s-\sigma_1|^2 |\rho-s|} \leq
\frac{(|\rho- s|+|s-\sigma_1|)^2}{|s-\sigma_1|^2 |\rho-s|} \leq
\frac{1}{|\rho- s|}+ \frac{2}{|s-\sigma_1|}+ \frac{|\rho-s|}{|s-\sigma_1|^2}\leq 1+\frac3a\, .
 $$
\end{itemize}
This proves that (\ref{eqn:cota}) is uniformly bounded.

\medskip

\subsection{Appendix 2: The exponent $d$ in Proposition \ref{prop:estimate} is best possible.}

\medskip
 We construct an example that has the sharp exponent. 

We construct 
a meromorphic function with convergence exponent $d=1$. More precisely, let $f$ be an entire function  
with zeros at $\rho=n^2 2^n i$, $n\geq 1$,
and with multiplicities $n_\rho=2^n$. Then $\sum n_\rho |\rho|^{-1}<\infty$.
The logarithmic derivative of such function is given by 
 $$
 g=\frac{f'}{f}=\sum \frac{n_\rho}{s-\rho}
 $$

Now let us see that it is not controlled as $|f'/f|\leq C |s|^{1-\epsilon}$ with $\epsilon>0$.
For this take $s= c+ k^2 2^k i$, $k$ a fixed integer, $c>0$. We decompose
 $$
 g(s)=\sum_{n=1}^{k-1} \frac{2^n}{c+ (k^2 2^k -n^22^n)i} +\frac{2^k}{c} + \
\sum_{n=k+1}^{\infty} \frac{2^n}{c+ (k^2 2^k -n^22^n)i} 
 $$
The first term is bounded by
 $$
 \sum_{n=1}^{k-1} \frac{2^n}{k^2 2^k-(k-1)^2 2^{k-1}} \leq
\frac{2^{k-1}}{2^{k-1}(k^2+2k-1)} < C_0,
 $$
for some universal constant.
The third term is bounded by 
$$ 
\sum \frac{2^n}{n^22^n-k^2 2^k}  \leq \sum \frac{2^n}{n^2 2^{n-1}} < C_1,
$$
for another universal constant. Hence
$|g(s)|\geq \frac{2^k}{c}-C_0-C_1$. For fixed $c$, take $k$ large
enough. Then 
 $$
 \frac{|g(s)|}{|s|^{1-\epsilon}} \geq
\frac{2^k/c-C_0-C_1}{(c^2+ k^4 2^{k+1})^{(1-\epsilon)/2}
} \approx \frac{2^{\epsilon}}{c \, k^{2-2\epsilon}},
 $$
which gets as large as we wish.

\subsection{Appendix 3: Proof of Lemma \ref{lem:ccc}.}  

Fix $c>\sigma_1$ and let $m_0$ be the minimal integer such that
$$
\left | (c+it)^{-m_0}  \frac{f'}{f} (c+it) \right |\in L^1(\RR ) \ . 
$$ 
Consider the holomorphic function 
 $$
g(s)=s^{-m_0} \frac{f'(s)}{f(s)}
$$
on the right half-plane $\Re s\geq c$.
The function $F(t)=g(c+it)$ satisfies the conditions of the Representation Theorem 6.5.4 in \cite{B} with 
$\alpha =0$, $c=0$, and we get using the last inequality of that Theorem
$$
\log |g(c'+iu)| \leq (c'-c)\pi^{-1} \int_\RR \frac {\log |g(c+it)| }{(t-u)^2+(c'-c)^2} \ dt \ .
$$ 
Now taking the exponential and using Jensen's convexity inequality we get
$$
|g(c'+iu)|\leq (c'-c)\pi^{-1} \int_\RR \frac {|g(c+it)| }{(t-u)^2+(c'-c)^2} \ dt \ .
$$
Now Fubini gives 
$$
 \int_\RR |g(c'+iu)|du \leq (c'-c)\pi^{-1} \int_\RR \left(
 |g(c+it)| \int_\RR\frac{1}{(t-u)^2+(c'-c)^2} \ du \right) dt =
 \int_\RR |g(c+it)| dt <\infty.
 $$

\end{document}